%% file: gsm.tex
\documentclass[11pt,a4paper, DIV15, oneside]{scrartcl}
\pdfoutput=1
\usepackage[english]{babel}
\usepackage[utf8x]{inputenc}

\usepackage[T1]{fontenc}
\usepackage[active]{srcltx}

\RequirePackage{amsthm,amsmath,natbib}
\RequirePackage{amssymb,amstext}
\RequirePackage{hypernat}
\usepackage[pdfpagemode=UseOutlines ,plainpages=false
,hypertexnames=false ,pdfpagelabels
,hyperindex=true,colorlinks=true]{hyperref}
	\makeatletter%
	\Hy@breaklinkstrue%
	\makeatother%
\usepackage{color}
\definecolor{darkred}{rgb}{0.6,0.0,0.1}
\definecolor{darkgreen}{rgb}{0,0.5,0}
\definecolor{darkblue}{rgb}{0,0,0.5}
\hypersetup{colorlinks ,linkcolor=darkblue ,filecolor=darkgreen
,urlcolor=darkblue ,citecolor=black}

\usepackage{abbrev-package}
\usepackage{abbrev-gsm}

\newtheoremstyle{mysc}
  {3pt}
  {3pt}
  {\it}
  {}
  {\sffamily\bfseries}
  {}
  {.5em}
  {}

\newtheoremstyle{myex}
  {3pt}
  {3pt}
  {\rm}
  {}
  {\sffamily\bfseries}
  {}
  {.5em}
  {}

\theoremstyle{mysc}\newtheorem{prop}{Proposition}[section]
\theoremstyle{mysc}\newtheorem{assumption}[prop]{Assumption}
\theoremstyle{mysc}\newtheorem{coro}[prop]{Corollary}
\theoremstyle{mysc}\newtheorem{theo}[prop]{Theorem}
\theoremstyle{mysc}\newtheorem{defin}[prop]{Definition}
\theoremstyle{mysc}\newtheorem{lem}[prop]{Lemma}
\theoremstyle{mysc}
\theoremstyle{myex}
\theoremstyle{myex}
\theoremstyle{myex}\newtheorem{illu}[prop]{Illustration}
\theoremstyle{myex}

\numberwithin{equation}{section}

\makeatletter%
\def\@fnsymbol#1{\ensuremath{\ifcase#1\or * \or \star \or 1 \or 2\or 3\or  , \or
g\or h\or i\else\@ctrerr\fi}}%
\makeatother%

\author{{\sc Jan Johannes}
  \and {\sc Maik
    Schwarz}\thanks{Corresponding author. Institut de
    statistique, biostatistique et sciences actuarielles, Voie du
    Roman Pays 20, bo\^ite L1.04.01, B-1348 Louvain-la-Neuve,
    Belgique, e-mail: \url{maik.schwarz@uclouvain.be}}
}

\title{Adaptive Gaussian inverse regression\\ with partially unknown operator}

\begin{document}
\date{Université catholique de Louvain}
\dedication{\today}
\maketitle

\begin{abstract}

\input{gsm_abstract.tex}
\end{abstract}

{\footnotesize

\begin{tabbing}
\noindent\emph{AMS 2010 subject classifications:} \= 62G05,
62G08.\\[.2ex]
\end{tabbing}
\begin{tabbing}
\noindent \emph{Keywords:} \=  Gaussian sequence space model, minimax
theory, adaptive nonparametric estimation, \\
\> model selection, Lepski's
method, Sobolev spaces, mildly and severely ill-posed inverse problems
\end{tabbing}
{\sf\textbf{Acknowledgments}} This work was supported by the IAP research network no.\ P6/03 the
Belgian Government (Belgian Science Policy) and by the ``Fonds
Sp\'eciaux de Recherche'' from the Universit\'e catholique de
Louvain. 
}


\input{gsm_intro.tex}

\input{gsm_minimax.tex}

\input{gsm_adaptive.tex}

\appendix

\section{Proofs}
\label{sec:proofs}

\input{gsm_minimax_proofs.tex}

\input{gsm_adaptive_proofs.tex}

\input{gsm_auxiliary.tex}

\bibliographystyle{apalike}
\bibliography{dr}

\end{document}

%% file: gsm_abstract.tex
\noindent This work deals with the ill-posed inverse problem of reconstructing a
function $\sol$ given implicitly as the solution of $\im = \Op\sol$,
where $\Op$ is a compact linear operator with unknown singular values
and known eigenfunctions. We observe the function $\im$ and the
singular values of the operator subject to Gaussian white noise with
respective noise levels $\nY$ and~$\nX$.

\noindent  We develop a minimax theory in terms of both noise levels and propose
 an orthogonal series estimator attaining the minimax rates. This
 estimator requires the optimal choice of a dimension parameter
 depending on certain characteristics of $\sol$ and $\Op$.  This work
 addresses the fully data-driven choice of the dimension parameter
 combining model selection with Lepski's method.  We show that the
 fully data-driven estimator preserves minimax optimality over a wide
 range of classes for $\sol$ and $\Op$ and noise levels $\nY$ and
 $\nX$. The results are illustrated considering Sobolev spaces and
 mildly and severely ill-posed inverse problems.


%% file: gsm_intro.tex
\section{Introduction}
\label{sec:introduction}

Let $(H,\skalar_H)$ and $(G,\skalar_G)$ be separable Hilbert spaces
and $\Op$ a compact linear operator from $H$ to~$G$ with unknown
singular values. This work deals with the reconstruction of a function
$\sol\in H$ given noisy observations of the image $\im=\Op\sol$ on the
one hand and of the unknown sequence of singular values
$\ev=(\ev_j)_{j\in\N}$ on the other hand. In other words, we consider a
statistical inverse problem with partially unknown operator.  There is
a vast literature on statistical inverse problems. For the case where
the operator is fully known, the reader may refer to \cite{JS:90}, 
\cite{MairRuymgaart1996}, \cite{MP:01}, and \cite{CGP:02}  and the
references therein. A typical illustration of such a situation is a
deconvolution problem (cf.\  \cite{0729.62089}, 
\cite{SC:90}, and \cite{Fan:91} among many others). For a more detailed discussion and
motivation of the case of a partially unknown operator which we
consider in this work, we refer the reader to \cite{CH:05}.
\cite{Efr:97} and \cite{Neu:97}  consider such a setting in the
particular context of a density deconvolution problem.

Let us describe in more detail the model we are going to consider. We suppose that $\Op$ admits
a singular value decomposition $(\ev_j,\bH_j,\bG_j)_{j\in\N}$ as
follows. Denote by $\Op^*$ the adjoint operator of $\Op$. Then,
$\Op^*\Op$ is a compact operator on $H$ with eigenvalues
$(\ev_j^2)_{j\in\N}$ whose associated orthonormal basis of
eigenfunctions $\{\bH_j\}$ we suppose to be known. Analogously, the
operator $\Op\Op^*$ has eigenvalues $(\ev_j^2)_{j\in\N}$ and known
orthonormal eigenfunctions $\bG_j = \normV{\Op\bH_j}_G^{-1}\Op\bH_j$
in $G$.  Projecting the inverse problem $\im=\Op\sol$ on the
eigenfunctions, we obtain the system of equations
$\fou{\im}_j:=\skalarV{\im,\bG_j}_G = \ev_j \skalarV{\sol,\bH_j}_H$
for $j\in\N$. As the operator $\Op$ is compact, the sequence of
singular values tends to zero and the inverse problem is called
ill-posed.

The solution $\sol$ is characterized by its coefficients $\fou{\sol}_j
:= \skalarV{\sol,\bH_j}_H$. Our objective is their estimation based on
the following observations:
  \begin{equation}
    \label{eq:model}
    Y_j = \fou{\im}_j
          + \sqrt{\nY}\,\errY_j
        = \ev_j\fou{\sol}_j 
          + \sqrt{\nY}\errY_j \qquad \text{and}\qquad
    X_j = \ev_j + \sqrt{\nX}\,\errX_j\qquad (j\in\N),
  \end{equation}
  where the $\errY_j,\errX_j$ are iid.\ standard normally distributed
  random variables and $\nY,\nX\in (0,1)$ are noise levels. Thus we
  represent the problem at hand as a hierarchical Gaussian sequence
  space model.

 Of course $\sol$ can only be reconstructed from such
  observations if all the $\ev_j$ are non-zero which is the case if
  and only if the operator $\Op$ is injective.  We assume this from
  now on, which allows us to write $\sol = \sum_{j=1}^\infty \fou{\im}_j\ev_j^{-1}\,\bH_j$.
Hence, an orthogonal
series estimator of $\sol$ is a natural approach:
\begin{equation*}\label{eq:12}
 \qquad\qquad \hsol_k := \sum_{j=1}^k
 \frac{Y_j}{X_j}\1_{[X_j^2\geq\,\sigma]}\, \bH_j.
\end{equation*}
The threshold using the indicator function
accounts for the uncertainty caused by estimating the $\ev_j$ by
$X_j$. It corresponds to $X_j$'s noise level as an estimator of
$\ev_j$, which is a natural choice \citep[cf.][p.310f.]{Neu:97}. Note
that $\hsol_k$ depends on a dimension parameter $k$ whose choice
essentially determines the estimation accuracy. Its optimal choice
generally depends on both unknown sequences $(\fou{\sol}_j)$
and~$(\ev_j)$.  Our purpose is to establish an adaptive estimation
procedure for the function $\sol$ which does not depend on these
sequences.  More precisely, assuming that the solution and the
operator belong to given classes $\sol\in\solclassletter$ and
$\Op\in\Opclassletter$, respectively, we shall measure the accuracy of
an estimator $\tsol$ of $\sol$ by the maximal weighted risk
$\cR_\hw(\tsol,\solclassletter,\Opclassletter) :=
\sup_{\sol\in\solclassletter}\sup_{\Op\in\Opclassletter}\Ex\normV{\tsol-\sol}_\hw^2$
defined \wrt some weighted norm $\norm_\hw :=
\sum_{j\in\N}\hw_j|\fou{\cdot}_j|^2$, where $\hw := (\hw_j)_{j\in\N}$
is a strictly positive weight sequences. This allows us to quantify
the estimation accuracy in terms of the mean integrated square error
(MISE) not only of $\sol$ itself, but as well of its derivatives, for
example. Given observations $Y=(Y_j)_{j\in\N}$ and $X =
(X_j)_{j\in\N}$ with respective noise levels $\nY$ and $\nX$ according
to~\eqref{eq:model}, the minimax risk \wrt the classes
$\solclassletter$ and $\Opclassletter$ is then defined as
$\cR^*_\hw(\nY,\nX,\solclassletter,\Opclassletter) :=
\inf_{\tsol}\cR_\hw(\tsol,\solclassletter,\Opclassletter)$, where the
infimum is taken over all possible estimators $\tsol$ of $\sol$. An
estimator $\hsol$ is said to attain the minimax rate or to be minimax
optimal \wrt $\solclassletter$ and $\Opclassletter$ if there is a
constant $C>0$ depending on the classes only such that
$\cR_\hw(\hsol,\solclassletter,\Opclassletter)\leq C\,
\cR^*_\hw(\nY,\nX,\solclassletter,\Opclassletter)$ for all
$\nY,\nX\in(0,1)$. An estimation procedure which is fully data-driven
and minimax optimal for a wide range of classes $\solclassletter$ and
$\Opclassletter$ is called \textit{adaptive}.

\medskip

In the next section, we show that for a wide range of classes
$\solclassletter$ and $\Opclassletter$ the orthogonal series
estimator~$\hsol_\kstar$ attains the minimax rate for an optimal
choice $\kstar$ of the dimension parameter. We illustrate this result
considering subsets of Sobolev spaces for~$\solclassletter$ and
distinguishing two types of operator classes~$\Opclassletter$
specifying the decay of the singular values: If~$(\ev_j)$ decays
polynomially, the inverse problem is called mildly ill-posed and
severely ill-posed if they decay exponentially.  However,~$\kstar$ is
chosen subject to a classical variance-squared-bias trade-off and
depends on properties of both classes~$\solclassletter$
and~$\Opclassletter$ which are unknown in general.

\medskip The last section is devoted to the development of a
data-driven choice $\hk$ of $k$, following the general model selection
scheme \citep[cf.]{BBM:99}. This methodology requires the careful
choice of a contrast function and a penalty term. In this work, we
will use a contrast function inspired by the work of \cite{GL:10} who
consider bandwidth selection for kernel estimators.  Given a random
sequence $(\hpen_k)_{k\geq 1}$ of penalties, a random set
$\{1,\dots,\hK_{\nY,\nX} \}$ of admissible dimension parameters and the random
sequence of contrasts
\begin{equation}\label{eq:con}
\hcon_k := \max_{k\leq j \leq \hK_{\nY,\nX}} \Big\{\normV{\hsol_j -
  \hsol_k}^2_\hw - \hpen_j   \Big\}\qquad (k\in\N).
\end{equation}
The dimension parameter is selected as the 
minimizer\footnote{For a sequence $(b_k)_{k\in\N}$ attaining a minimal
  value on $N\subset\N$, let $\argmin\limits_{n\in N}b_n := \min\{ n\in N\;|\;b_n \leq
  b_k \; \forall k\in N\}$.}
 of a penalized contrast
 \begin{equation}\label{eq:khat}
\hk := \argmin_{1\leq k \leq \hK}\Big\{ \hcon_k + \hpen_k    \Big\}.
\end{equation}
We assess the accuracy of the fully data-driven estimator $\hsol_\hk$
deriving an upper bound for
$\cR_\hw(\hsol_\hk,\solclassletter,\Opclassletter)$.  Obviously this
upper bound heavily depends the random sequence $(\hpen_k)$ and the
random upper bound $\hK$.  However, we construct these objects in such
a way that the resulting fully data-driven estimator $\hsol_\hk$ is
minimax optimal over a wide range of classes and thus adaptive. 
The more technical proofs and some auxiliary results are deferred to
the appendix.

\cite{HR:08} also study adaptive estimation in linear inverse
problems, but their method is limited to mildly ill-posed inverse
problems with known degree of ill-posedness. Also, the theoretical
framework is quite different: they focus on sparse representations and
therefore consider estimators based on wavelet thresholding and show
their rate-optimality and adaptivity properties over Besov spaces with
respect to the corresponding norms.

Adaptive estimation in a hierarchical Gaussian sequence space model
has previously been considered by \cite{CH:05}. Though, the authors
restrict their investigation to the mildly ill-posed case and to noise
levels satisfying $\nX\leq \nY$. The new approach presented in this
paper has the advantage of not requiring such restrictions.  On the
contrary, the influence of the two noise levels on the estimation
accuracy is characterized. Moreover, the estimator presented in this
paper can attain optimal convergence rates independently of whether
the underlying inverse problem is mildly or severely ill-posed, for
example, even when~$\nY\ll\nX$. This is an
important feature in applications where the reduction of the noise
level $\nX$ can be costly. In (satellite or medical) imaging, for
example, the observation of the sequence $X$ may correspond to
calibration measurements. In order to achieve an adequately high
precision of these measures as to reduce the noise level $\nX$
sufficiently, one might have to repeat expensive experiments. It is
thus desirable to know how the estimator performs when~$\nX$ exceeds
$\nY$.


%% file: gsm_minimax.tex
\section{Minimax}
\label{sec:minimax}

In this section we develop a minimax theory for Gaussian inverse
regression with respect to the classes
\begin{align*}\label{eq:classes}
  \begin{split}
    \solclass &:=\bigg\{ h \in H\;\Big|\;\sum_{j\in\N} \solw_j
    |[h]_j|^2
    =:\normV{h}_\solw^2\leq \solr\bigg\}\mbox{ and } \\
    \Opclass &:=\bigg\{ T\in C(H,G) \;\Big|\; \text{The eigenvalues
      $\{u_j\}$ of $T^*T$ satisfy}\; 1/\Opr \leq
    \frac{u_j^2}{\Opw_j}\leq \Opr\quad \forall\, j\in\N\bigg\},
  \end{split}
\end{align*}
where $C(H,G)$ denotes the set of all compact linear
operators from $H$ to $G$ having $\{\bH_j\}$ and $\{\bG_j\}$ as eigenfunctions, respectively. The minimal regularity conditions on the
solution, the operator and the weighted norm $\norm_\hw$ which we need
in this section are summarized in the following assumption. 
\begin{assumption}\label{ass:minreg} Let $\solw:=(\solw_j)_{j\in \N}$,
  $\hw:=(\hw_j)_{j\in \N}$ and $\Opw:=(\Opw_j)_{j\in \N}$ be strictly
  positive sequences of weights with
  $\solw_1=\hw_1=\Opw_1 =1$ such that  $\hw/\solw$ and $\Opw$ are
  non-increasing, respectively.
\end{assumption}

\begin{illu}\label{illu:ooo}
As an illustration of the results below, we will consider weight sequences
$\solw_j = j^{2p}$, for which $\solclass$ is a Sobolev space of
$p$-times differentiable functions if we
consider the trigonometric basis in $H=L^2[0,1]$. As for the operator, we will
distinguish the cases $\Opw_j= j^{-2b}$, referred to as 
\textit{mildly ill-posed} (\textbf{[m]}) and
$\Opw_j = \exp(-j^{2b})$, the \textit{severely ill-posed} case
(\textbf{[s]}). Concerning the weighted norm, we will consider sequences\footnote{$b_\rho\sim c_\rho$ means
    that $\lim_{\rho\to 0} b_\rho/c_\rho$ exists in $(0,\infty)$.  }
$\hw_j\sim j^{2s}$, such that $\normV{f}_\hw = \normV{f^{(s)}}_{L^2}$ for all
$f\in\solclass$. We will assume that $b\geq 0$ and $p\geq s \geq 0$,
such that Assumption~\ref{ass:minreg} is satisfied. 
\end{illu}
The following result states lower risk bounds for the estimation of
$\sol$ and thus describes the complexity of the problem. 
\begin{theo} \label{theo:lower}  Suppose that we observe
sequences  $Y$ and $X$ according to the model~\eqref{eq:model}. Consider sequences $\hw$, $\solw$, and $\Opw$ satisfying
  Assumption~\ref{ass:minreg}.
  For all $\nY,\nX\in
  (0,1)$, define
  \begin{equation}
\label{def:psin}
   \rho_{k,\nY} := \max\Bigl(
      \frac{\hw_k}{\solw_k},
      \sum_{j=1}^k\frac{\nY\hw_j}{\Opw_j}\Bigr),\quad
      \krate_\nY:=\min_{k\in\N} \rho_{k,\nY}, \quad
      \kstar:=\argmin\limits_{k\in\N} \rho_{k,\nY},\quad
      \urate_\nX:=\max_{k\in\N} \Bigr\{
      \frac{\hw_k}{\solw_k}\min\Bigl(1,\frac{\nX}{\Opw_k}\Bigr)\Bigr\}.
    \end{equation}
  If   $\eta:=\inf_{n\in\N}\{\krate_\nY^{-1}\min(\hw_{\kstar}\solw_{\kstar}^{-1},
  \sum_{l=1}^{\kstar}{\nY\hw_l}{(\Opw_l)}^{-1}) \}>0$, then
  \begin{align*}
\inf_{\tsol}    \cR_\hw(\tsol,\solclass,\Opclass)
    \geq
\frac{1}{4d} 
\min(\eta,r)  \min(r,1/(2d), (1-d^{-1/2})^2)\;
    \max(\krate_\nY, \urate_\nX),
  \end{align*}
where the infimum is to be taken over all possible estimators $\tsol$
of $\sol$. 
\end{theo} It is noteworthy that apart from the unwieldy constant, the
lower bound is given by two terms ($\krate_\nY$ and $\urate_\nX$),
each of which depending only on one noise level. We show in the proof
that $\krate_\nY$ is actually, up to a constant, a lower risk bound
uniformly for any known operator $\Op$ in the class $\Opclass$. Hence,
in this case no supremum over the class $\Opclass$ would be needed.
The term $\urate_\nX$ only arises if the operator is unknown in
$\Opclass$. The proof of this lower bound is based on a comparison of
different inverse problems with different operators in $\Opclass$,
whence the supremum over this class.  The term $\urate_\nX$ quantifies
to which extent the additional difficulty arising from the preliminary
estimation of the eigenvalues $\ev_j$ influences the possible
estimation accuracy for $\sol$: As long as $\krate_\nY\geq
\urate_\nX$, the same lower bound as in the case of known eigenvalues
holds. Otherwise, the lower bound increases.  Notice further that the
term $\rho_{k,\nY}$ above corresponds to the MISE of the orthogonal
series estimator $\hsol_k$ in the case of known eigenvalues $\ev_j$,
and $\kstar$ is its minimizer \wrt $k$. Under classical smoothness
assumptions, the rates and~$\kstar$ take the following forms.

\begin{illu}\label{illu:chika}
  In the special cases defined in Illustration~\ref{illu:ooo} above,
  the rates from~\eqref{def:psin} are\\[1ex]
\textbf{[m]}\qquad $ \krate_\nY \sim \nY^{2(p-s) / (2p+2b+1)},\qquad
     \kstar \sim \nY^{-1 / (2p+2b+1)}, \qquad
\urate_\nX \sim \nX^{((p-s)\wedge b)/b} $   \\[1ex]
\textbf{[s]} \qquad $
\krate_\nY\sim |\log\nY|^{(p-s)/b},\qquad
     \kstar\sim |\log\nY|^{1/(2b)},\qquad
   \urate_\nX\sim |\log\nX|^{-(p-s)/b}$.
\end{illu}

\vspace{1em}
 The following theorem shows that the orthogonal series estimator
  $\hsol_\kstar$ with optimal parameter $\kstar$ given in~\eqref{def:psin}
  actually attains the lower risk bound up to a constant and is thus
  minimax optimal.
\begin{theo}\label{theo:upper}Under the assumptions of
  Theorem~\ref{theo:lower}, the estimator $\hsol_\kstar$ satisfies for all $\nY,\nX \in (0,1)$ 
  \begin{align*}
    \sup_{\sol\in\solclass}\sup_{\Op\in\Opclass} \left\{
      \Ex\normV{\hsol_{\kstar}-\sol}^2_\hw\right\}
\leq 4(6\Opr + \solr)\,\max(\krate_\nY,\urate_\nX).
  \end{align*}
\end{theo} 

\vspace{1em}
 To conclude this section, let us summarize the resulting optimal convergence rates
   under the classical smoothness
  assumptions introduced in Illustration~\ref{illu:ooo}. In order to
characterize the influence of the second noise level $\nX$, we
consider it as a function of the first noise level $\nY$.
\begin{illu}\label{illu:rates_minimax}
Let $(\nX_\nY)_{\nY\in(0,1)}$ be a noise level in $X$ depending on the
noise level $\nY$ in $Y$. \\[1ex]
\textbf{[m]} Let $p> 1/2$, $b>1$, and $0\leq s\leq p$. 
If $q_1:=\lim\limits_{\nY\to 0} \nY^{-2((p-s)\vee b)/(2p+2b+2)}\nX_\nY$ 
exists\footnote{The limit <<$\infty$>> meaning strict divergence is
  authorized.}, then 
\[\sup_{\sol\in\solclass}\sup_{\Op\in\Opclass}\E\normV{\hsol^{(s)}_\kstar - \sol^{(s)}}^2_{L^2} = 
\begin{cases}
  O(\nY^{2(p-s)/(2p+2b+1)})   &   \text{if $q_1<\infty$} \\
  O(\nX_\nY^{((p-s)\wedge b)/b})   &   \text{otherwise.}
\end{cases} \]
\textbf{[s]} Let $p>1/2$,$b>0$ and $0\leq s \leq p$. 
If $q_2:=\lim\limits_{\nY\to 0} |\log\nY| \, |\log\nX_\nY|^{-1}$
exists, then
\[\sup_{\sol\in\solclass}\sup_{\Op\in\Opclass}\E\normV{\hsol^{(s)}_\kstar - \sol^{(s)}}^2_{L^2} = 
\begin{cases}
  O(|\log\nY|^{(p-s)/b})   &   \text{if $q_2<\infty$} \\
  O(|\log\nX_\nY|^{(p-s)/b})   &   \text{otherwise.}
\end{cases}\]
\end{illu}
This illustration shows that often the same optimal rates as in the
case of known eigenvalues hold even when $\nY<\nX$.


%% file: gsm_adaptive.tex
\section{Adaptation}
\label{sec:adaptive}

In this section, we construct a fully data-driven estimator of
$\sol$ following the procedure sketched in~\eqref{eq:con}
and~\eqref{eq:khat}. The following Lemma will be our key tool when
controlling the risk of the adaptive estimator.

\begin{lem}\label{lem:keylemma}
  Let $\pen$ be an arbitrary positive sequence and
  $K\in\N$. Consider the sequence $\con$  of contrasts $\con_k := \max_{k\leq j \leq K} \Big\{\normV{\hsol_j -
  \hsol_k}^2_\hw - \pen_j   \Big\}$ and $\tk:= \argmin_{1\leq j \leq
  K}\{\con_j+\pen_j\}$. Let further $(t)_+ := (t\vee 0)$. If~$(\pen_1,\dots,\pen_K)$ is non-decreasing,
then we have for all $1\leq k \leq K$ that
\begin{equation}
\normV{\hsol_\tk-\sol}^2_\hw\leq 7\pen_k + 78\bias^2_k + 
 42\max_{1\leq j\leq K}\Big(\normV{\hsol_j - \sol_j}^2_\hw -
 \frac{1}{6}\pen_j   \Big)_+,\label{eq:keylem}
\end{equation}
where we denote by  $\sol_j := \sum_{k=1}^j \fou{\sol}_k\,\bH_k$ the
projection of $\sol$ on the first $j$ basis vectors in $H$ and
by
$\bias_k:= \sup_{j\geq k}\normV{\sol - \sol_j}_\hw$ the bias due to the
projection.
\end{lem}
\proof
In view of the definition of $\tk$, we have for all $1\leq k\leq K$
that
\begin{align}\label{eq:key:1}
  \begin{split}
    \normV{\hsol_\tk - \sol}_\hw^2 &\leq 3 \Big\{\normV{\hsol_\tk -
      \hsol_{k\wedge\tk}}_\hw^2 + \normV{\hsol_{k\wedge\tk} -
      \hsol_{k}}_\hw^2
    +  \normV{\hsol_k - \sol}_\hw^2   \Big\} \\
    &\leq 3\Big\{ \con_k + \pen_\tk + \con_\tk + \pen_k
    +\normV{\hsol_k -
      \sol}_\hw^2      \Big\}\\
    &\leq 6\Big\{ \con_k + \pen_k \Big\} +3\normV{\hsol_k -
      \sol}_\hw^2
  \end{split}
\end{align}
Since $(\pen_1,\dots,\pen_K)$ is non-decreasing and
$4\bias_k^2\geq\max_{k\leq j \leq K}\normV{\sol_k-\sol_j}_\hw^2$, we have
\[\con_k\leq 6 \max_{1\leq j \leq K}\Big(\normV{\hsol_j -
  \sol_j}_\hw^2 - \frac{1}{6} \pen_j    \Big)_+ + 12 \bias_k^2.\]
It easily verified that for all $1\leq k \leq K$ we have
\[\normV{\hsol_k - \sol}_\hw^2 \leq \frac{1}{3}\pen_k + 2\bias_k^2 + 2\max_{1\leq j \leq K}\Big(\normV{\hsol_j -
  \sol_j}_\hw^2 - \frac{1}{6} \pen_j    \Big)_+.\]
The result follows combining the last estimates with~\eqref{eq:key:1}.
\qed

The Lemma being valid for any upper bound $K$ and any monotonic
sequence of penalties $\pen$, we need to specify our choice.  Let us first
define some auxiliary quantities required in the construction of the
random penalty sequence $\hpen$ and the upper bound $\hK$.

\begin{defin}\label{def:deltagen}
For any sequence $\alpha := (\alpha_j)_{j\in\N}$, define
\ii
\begin{enumerate}
\item $\Delta_k^\alpha := \max_{1\leq j \leq k} \hw_j\,\alpha_j^{-2} $
  \qquad and \qquad $\delta_k^\alpha := k\Delta_k^\alpha
  \frac{\log(\Delta_k^\alpha \vee (k+2))}{\log(k+2)}$;
\item given $\hw_k^+ := \max_{1\leq j\leq k}\hw_j$, $N^\circ_\nY
  := \max\{1\leq N \leq \nY^{-1} \;|\; \hw_N^+ \leq \nY^{-1}\}$, \\and 
$v_\nX:= (8\log(\log(\nX^{-1}+20)))^{-1}$, let
\[N_\nY^\alpha := \min\Big\{ 2\leq j \leq N^\circ_\nY \;\Big|\;
     \frac{\alpha_j^2}{j\hw_j^+}  \leq \nY|\log\nY|  \Big\} -1
\quad\text{and} \quad
  M_\nX^\alpha :=  \min\Big\{ 2\leq j \leq \nX^{-1} \;\Big|\;
     \alpha_j^2 \leq \nX^{1-v_\nX}  \Big\} -1,\]
and $K_{\nY,\nX}^\alpha := N_\nY^\alpha \wedge M_\nX^\alpha$. If the
defining set is empty, set $N_\nY^\alpha = N^\circ_\nY$ or $M_\nX^\alpha = \lfloor\nX^{-1}\rfloor$, respectively.
\end{enumerate}
\end{defin}

Choosing appropriate sequences $\alpha$, these quantities allow us
define the random penalty term needed for the data-driven choice of
$k$ as well as its deterministic counterpart which will be used in the
control of the risk.

Using this definition and denoting by $X$ the sequence of random
variables $(X_j)_{j\in\N}$, define 
  \begin{equation}\label{Khatpenhat}
\hK_{\nY,\nX}:=
  K^X_{\nY,\nX} \quad\text{and}\quad \hpen_k:=600\delta_k^X\,\nY.
\end{equation}
Substituting these definitions in~\eqref{eq:con} and~\eqref{eq:khat}
yields a choice of the dimension parameter $k$ depending exclusively
on the observations and the noise levels, but not on any underlying
smoothness classes. 

Consider the upper risk bound in Lemma~\ref{lem:keylemma}. In order to
control the risk of the data-driven estimator, we decompose it with
respect to an event on which the randomized quantities $\hpen_k$ and
$\hK_{\nY,\nX}$ are close to some deterministic counterparts
$\pen^a_k$, $K^-_{\nY,\nX}$, and $K^+_{\nY,\nX}$ to be defined below
in Propositions~\ref{prop:concentration}
and~\ref{prop:schachtel}. More precisely, consider the event
\[\mho_{\nY,\nX} := \{\pen^a_k\leq \hpen_k \leq 30\pen^a_k \quad
\forall\;1\leq k \leq K^+_{\nY,\nX}\} \cap \{K^-_{\nY,\nX} \leq
\hK_{\nY,\nX} \leq K^+_{\nY,\nX}\}\] 
and the corresponding risk decomposition
\begin{equation}
\E\normV{\hsol_\hk-\sol}^2_\hw  =
\E\normV{\hsol_\hk-\sol}^2_\hw\1_{\mho_{\nY,\nX}}
+\E\normV{\hsol_\hk-\sol}^2_\hw \1_{\mho_{\nY,\nX}^c}.\label{eq:risk_mho}
\end{equation}

As the random sequence $\hpen_k$ is non-decreasing in $k$ by construction, we may
apply Lemma~\ref{lem:keylemma} and obtain for every $1\leq k \leq \hK_{\nY,\nX}$
\[\normV{\hsol_\hk-\sol}^2_\hw\leq 7\,\hpen_k + 78\bias^2_k + 
 42\max_{1\leq j\leq \hK_{\nY,\nX}}\Big(\normV{\hsol_j - \sol_j}^2_\hw -
 \frac{1}{6}\hpen_j   \Big)_+.\]
On the event $\mho_{\nY,\nX}$, this implies that
\begin{equation}\label{eq:risk_on_mho}
\E\normV{\hsol_\hk-\sol}^2_\hw\1_{\mho_{\nY,\nX}}
\leq 420\min_{1\leq k \leq K^-_{\nY,\nX}}\{\max( \pen^a_k, \bias^2_k)\}
+ 42\max_{1\leq j\leq K^+_{\nY,\nX}}\E\Big(\normV{\hsol_j - \sol_j}^2_\hw -
 \frac{1}{6}\pen^a_j   \Big)_+.
\end{equation}

The second term in the last inequality is controlled uniformly over
$\solclass$ and $\Opclass$ by the following Proposition.

\begin{prop}\label{prop:concentration}
  Given $\Op\in\Opclass$ with singular values $a:=(a_j)_{j\in\N}$, let
 $\sqrt{4\Opr\Opw} := (\sqrt{4\Opr\Opw_j})_{j\in\N}$ and
 define $K^+_{\nY,\nX}:= K^{\sqrt{4\Opr\Opw}}_{\nY,\nX}$,
  $M^+_{\nY,\nX}:= M^{\sqrt{4\Opr\Opw}}_{\nY,\nX}$, and
  $\pen^a_k:=60\delta_k^\ev\,\nY$ using Definition~\ref{def:deltagen}.
  There is a constant $C>0$ depending only on the class $\Opclass$
  such that
\[\sup_{\sol\in \solclass}\sup_{\Op\in\Opclass}
 \E \Big[\max _{1\leq k \leq K^+_{\nY,\nX}}
 \Big(\normV{\hsol_k - \sol_k}^2_\hw - \frac{1}{6}\pen^a_k    \Big)_+
 \Big]  \leq C\,\Big\{\nY + \solr\urate_\nX + \nX\Big\}.\]
\end{prop}

Roughly speaking, the penalty term is an upper bound for the
estimator's variation. Typically, it can be chosen as a multiple of
the estimator's variance. Thus, inequality~\eqref{eq:keylem} actually
features a bias variance decomposition of the risk with an additional
third term which is controlled by the above proposition.
 
\begin{illu}\label{illu:pen}
   Note that for any operator $\Op\in\Opclass$ with sequence
  $(\ev_j)_{j\geq1}$ of singular
  values, the sequence $\delta^a$ appearing in the definition of the
  penalty term $\pen^a$ satisfies $(\Opr\zeta_\Opr)^{-1}\leq
  (\delta^\ev_j/\delta^\Opw_j)\leq\Opr\zeta_\Opr$ for all $j\in\N$, with  $\zeta_d =
 \log(3\Opr)/\log(3)$. In the special cases defined in Illustration~\ref{illu:ooo}
  above, the sequence $\delta^\Opw$ takes the following form:  \\[1ex]
  \textbf{[m]}\quad $\delta^\Opw_k\sim k^{2b+2s+1}$   \hspace{4em}
  \textbf{[s]}\quad $\delta^\Opw_k\sim k^{2b+2s+1}\exp(k^{2b})(\log
  k)^{-1}$.
\end{illu}
\bigskip

The next proposition ensures that the randomized upper bound and
penalty sequence behave similarly to their deterministic counterparts
with sufficiently high probability so as not to deteriorate the
estimation risk. In view of Proposition~\ref{prop:concentration}, this justifies the choice of the penalty.
\begin{prop}\label{prop:schachtel}
  Let $K^-_{\nY,\nX}:=
  K^{\sqrt{\Opw/(4\Opr)}}_{\nY,\nX}$ and $ M_\nX^+ :=
M_\nX^{\sqrt{4\Opr\Opw}}$ using Definition~\ref{def:deltagen} and  
 suppose that there is a constant $L>0$ depending only on $\Opw$ and $\Opr$ such that
 \begin{equation}
   \nX^{-7}\Opw_{\Ce^+_\nX+1}^{-1/2}\exp\left( - {\Opw_{\Ce^+_\nX+1}}/({72\,\nX\edfr})
    \right) \leq L \quad\text{for all}\quad \nX\in (0,1).\label{eq:Mplusone}
  \end{equation}
 Then, there is a constant $C>0$ 
depending only on the class $\Opclass$ such that
\[\sup_{\sol\in \solclass}\sup_{\Op\in\Opclass}\E[\normV{\hsol_\hk -
  \sol}^2_\hw \1_{\mho^c_{\nY,\nX}} ] \leq
C\,(1+\solr)\,\nX\quad\text{for all}\quad\nY,\nX\in(0,1).\] 
\end{prop}

Condition~\eqref{eq:Mplusone} is satisfied in particular under the classical
smoothness assumptions considered in the illustrations. 
We are finally prepared to state the upper risk bound of the
  fully data-driven estimator~$\hsol_\hk$ of~$\sol$, which is the main result of
  this article.

\begin{theo}\label{theo:adaptive} Under Assumption~\ref{ass:minreg} and supposing~\eqref{eq:Mplusone},
  there is a constant $C$ depending only on the class $\Opclass$ such that for all $\nY,\nX\in(0,1)$ the adaptive
  estimator $\hsol_\hk$ satisfies
\begin{equation*}
\cR_\hw(\hsol_\hk,\solclass,\Opclass)\leq C\,(1+\solr)\,\Big\{\min_{1\leq k \leq
  K_{\nY,\nX}^-}\{\max(\hw_k/\solw_{k},  \delta^\Opw_{k}\nY)\}
+\urate_\nX + \nY + \nX \Big\}.
\end{equation*}
\end{theo}
\proof Considering~\eqref{eq:risk_on_mho}, note that for all
$\Op\in\Opclass$, we have $\pen^\ev_k\leq
60\nY\Opr\zeta_\Opr\delta^\Opw_k$ with $\zeta_d =
\log(3\Opr)/\log(3)$. On the other hand, it is easily seen that for
all $\sol\in\solclass$, one has $\bias^2_k \leq
\solr\,(\hw_k/\solw_k)$. Thus, we can write
\[\sup_{\sol\in\solclass}\sup_{\Op\in\Opclass}\min_{1\leq k \leq
  K^-_{\nY,\nX}} \{\max(\pen^\ev_k,\bias^2_k)\} \leq C\,(1+\solr)\,
\min_{1\leq k \leq K^-_{\nY,\nX}}
\{\max(\hw_k/\solw_k,\delta^\Opw_k\nY) \} \] for some constant $C>0$
depending only on $\Opr$.  In view of~\eqref{eq:risk_mho}, the rest of
the proof is obvious using Propositions~\ref{prop:concentration}
and~\ref{prop:schachtel}.  \qed

A comparison with the lower bound from Theorem~\ref{theo:lower} shows
that this upper bound ensures minimax optimality of the adaptive
estimator $\hf_\hk$ only if 
\[\krate^\diamond_{\nY,\nX}:= \min_{1\leq k \leq K^-_{\nY,\nX}}
 \Big[\max\Big( \frac{\hw_k}{\solw_k} , \delta^\Opw_k\nY \Big)
 \Big]\]
is at most of the same order as $\max(\krate_{\nY},\urate_\nX)$, whence the following
corollary. 
\begin{coro}\label{cor:adapivoptimal}
  Under Assumption~\ref{ass:minreg} and if
  $\sup_{\nY,\nX\in(0,1)}\{\krate^\diamond_{\nY,\nX} /
  \max(\krate_\nY,\urate_\nX)\}<\infty$, we have 
\[ \cR_\hw(\hsol_\hk,\solclass,\Opclass) \leq C\,\cR^*_\hw(\solclass,\Opclass) \qquad \forall\; \nY,\nX\in(0,1).\]
\end{coro}
We conclude this article reconsidering the framework of the
  preceding Illustration~\ref{illu:rates_minimax}. Notice that the
  adaptive estimator is minimax optimal over a wide range of cases,
  even when $\nY<\nX$.
\begin{illu}\label{illu:adaptive}
Let $(\nX_\nY)_{\nY\in(0,1)}$ be a noise level in $X$ depending on the
noise level $\nY$ in $Y$ and suppose that the limits $q_1$ and
$q_2$ from Illustration~\ref{illu:rates_minimax} exist in the
respective cases. Some straightforward computations then show
  that the adaptive estimator attains the following rates of
  convergence.\\[1ex]
\textbf{[m]} If $p-s>b$, the adaptive estimator $\sol_\hk^{(s)}$
attains the optimal rates (cf.~Illustration~\ref{illu:rates_minimax}).
In case $p-s\leq b$, we have, supposing that
$q^v_1:=\lim\limits_{\nY\to
  0}\nY^{-2b/(2p+2b+1)}\nX_\nY^{1-v_{\nX_\nY}}$ exists,
 \[\sup_{\sol\in\solclass}\sup_{\Op\in\Opclass}\E\normV{\hsol^{(s)}_\hk
   - \sol^{(s)}}^2_{L^2} = 
 \begin{cases}
   O(\nY^{2(p-s)/(2p+2b+1)})    &  \text{if $q_1<\infty$ and $q^v_1<\infty$, }  \\
   O(\nX_\nY^{(p-s)/b}\nX_\nY^{-v_{\nX_\nY}}) & \text{otherwise.}
 \end{cases}
\]\\[1ex]
\textbf{[s]} The adaptive estimator attains the optimal rates. 
\end{illu}


%% file: gsm_minimax_proofs.tex
\subsection{Minimax theory (Section \ref{sec:minimax})}
\label{sec:proofs_minimax}

\subsubsection*{Lower risk bound}
\label{sec:lower-risk-bound}

\proofof{Theorem \ref{theo:lower}} 
The proof consists of two steps: (A) First, we show that $\krate_\nY$
yields a lower risk bound in the case where the eigenvalues $(\ev_j)$ of
the operator $\Op$ are known. (B) Then, we show that another lower
risk bound is given by $\urate_\nX$.

\textit{Step (A).} Given $\zeta:=\eta\min(\solr,
1/(2\Opr))$ and $\alpha_\nY:=\krate_\nY (\sum_{j=1}^{\kstar}
\nY\hw_j/\Opw_j)^{-1}$ we consider the function $\sol:= (\nY\zeta\alpha_\nY)^{1/2}\sum_{j=1}^{\kstar}\Opw_j^{-1/2}
\bH_j$. We are going to show that for any
$\theta:=(\theta_j)\in\{-1,1\}^{\kstar}$, the function $\sol_\theta:=
 \sum_{j=1}^{\kstar} \theta_j\fou{\sol}_j \bH_j$ belongs to
$\solclass$ and is hence a possible candidate for the solution. 

For a fixed $\theta$ and under the hypothesis that the solution is
$\sol_\theta$, the observation $Y_k$ is distributed according to
$\cN(\ev_k\fou{\sol_\theta}_k, \nY)$ for any $k\in\N$. 
We denote by $\PP_\theta$ the distribution of the resulting sequence
$\{Y_k\}$ and by $\E_\theta$ the expectation \wrt this
distribution.

Furthermore, for $1\leq j\leq\kstar$ and each
$\theta$, we introduce $\theta^{(j)}$ by $\theta^{(j)}_{l}=\theta_{l}$
for $j\neq l$ and $\theta^{(j)}_{j}=-\theta_{j}$.  The key argument of
this proof is the following reduction scheme. If $\tsol$ denotes an
estimator of $\sol$ then we conclude
\begin{align}\label{eq:11}
  \begin{split}
    \sup_{\sol\in\solclass} &\Ex\normV{\tsol -\sol}_\hw^2 \geq
    \sup_{\theta\in
      \{-1,1\}^{\kstar}} \Ex_\theta\normV{\tsol -\sol_\theta}_\hw^2\geq  \frac{1}{2^{{\kstar}}}\sum_{\theta\in \{-1,1\}^{2\kstar}}\Ex_\theta\normV{\tsol -\sol_\theta}_\hw^2\\
    &\geq \frac{1}{2^{{\kstar}}}\sum_{\theta\in \{-1,1\}^{\kstar}}\sum_{j=1}^{\kstar}\hw_j\Ex_{{\theta}}|[\tsol-\sol_\theta]_j|^2\\
    &= \frac{1}{2^{{\kstar}}}\sum_{\theta\in
      \{-1,1\}^{\kstar}}\sum_{j=1}^{\kstar}\frac{\hw_j}{2}\Bigl\{\Ex_{{\theta}}|[\tsol-\sol_\theta]_j|^2+\Ex_{{\theta^{(j)}}}|[\tsol-\sol_{\theta^{(j)}}]_j|^2
    \Bigr\}.
  \end{split}
\end{align}
Below we show furthermore that for all $\nY\in (0,1)$ we have
\begin{equation}\label{eq:443}
  \Bigl\{\Ex_{{\theta}}|[\tsol-\sol_\theta]_j|^2+\Ex_{{\theta^{(j)}}}|[\tsol-\sol_{\theta^{(j)}}]_j|^2
  \Bigr\} \geq 
  \frac{\nY\,\zeta\alpha_\nY}{2\Opw_j}.
\end{equation} 
Combining the last lower bound and the reduction scheme gives
\begin{equation*}
  \sup_{\sol\in\solclass} \Ex\normV{\tsol -\sol}^2_\hw \geq
  \frac{1}{2^{{\kstar}}}\sum_{\theta\in
    \{-1,1\}^{\kstar}}\sum_{j=1}^{\kstar}\frac{\hw_j}{2}
\frac{\nY\zeta\alpha_\nY}{2\Opw_j}
  = \frac{\zeta\alpha_\nY}{4}  \sum_{j=1}^{\kstar}
  \frac{\nY\hw_j}{\Opw_j}
  = \frac{\zeta\krate_\nY}{4},
\end{equation*}
which implies the lower bound given in the theorem by definition of $\zeta$.

To complete the proof, it remains to check \eqref{eq:443} and
$\sol_\theta\in\solclass$ for all $\theta\in\{-1,1\}^{\kstar}$. The
latter is easily verified if $\sol\in \solclass$, which can be seen
recalling that $\hw/\xdfw$ is non-increasing and noticing that the
definitions of $\zeta$, $\alpha_\nY$ and $\eta$ imply
$\normV{\sol}_\solw^2 \leq
\zeta\frac{\solw_{\kstar}}{\hw_{\kstar}}\alpha_\nY\Bigl(\sum_{j=1}^{\kstar}
\frac{\nY\hw_j}{\Opw_j}\Bigr) \leq { \zeta/\eta\leq \solr}$.
    
It remains to show~\eqref{eq:443}. Consider the Hellinger affinity
$\rho(\PP_1,\PP_{-1})= \int
\sqrt{d\PP_1\;d\PP_{-1}}$, then we obtain
for any estimator $\txdf$ of $\xdf$ that
\begin{align*}
  \rho({\PP_1,\PP_{-1}})&\leqslant \int
  \frac{|[\tsol-\sol_{\theta^{(j)}}]_j|}{|[\sol_\theta-\sol_{\theta^{(j)}}]_j|}\sqrt{d\PP_1\;d\PP_{-1}}
  + \int
  \frac{|[\tsol-\sol_{\theta}]_j|}{|[\sol_\theta-\sol_{\theta^{(j)}}]_j|}
  \sqrt{d\PP_1\;d\PP_{-1}}\\
  &\leqslant \Bigl( \int
  \frac{|[\tsol-\sol_{\theta^{(j)}}]_j|^2}{|[\sol_\theta-\sol_{\theta^{(j)}}]_j|^2}
  {d\PP_1}\Bigr)^{1/2} + \Bigl( \int
  \frac{|[\tsol-\sol_{\theta}]_j|^2}{|[\sol_\theta-\sol_{\theta^{(j)}}]_j|^2}
  d\PP_{–1}
\Bigr)^{1/2}.
\end{align*}
Rewriting the last estimate we obtain
\begin{equation}\label{pr:theo:lower:n:e2}
  \Bigl\{\Ex_{{\theta}}|[\tsol-\sol_\theta]_j|^2+\Ex_{{\theta^{(j)}}}|[\tsol-\sol_{\theta^{(j)}}]_j|^2
  \Bigr\} \geq  \frac{1}{2} |[\sol_\theta-\sol_{\theta^{(j)}}]_j|^2
  \rho^2(\PP_1,\PP_{-1}). \end{equation}
Next, we bound the Hellinger affinity
$\rho({\PP_1,\PP_{-1}})$ from below. 
Consider the Kullback-Leibler divergence of these two distributions first. 
The components of the two sequences corresponding to the
distributions $\PP_1$ and $\PP_{-1}$ are pairwise equally distributed
except for the $j$-th component. Thus, we have
$\log({d\PP_\theta}/{d\PP_\thetaj})=({2y_j\ev_j\theta_j\fou{\sol}_j}/{
\nY})$,
and taking the integral over $y_j$ \wrt $\PP_\theta$, we find
\[{KL}({\PP_1,\PP_{-1}}) = \frac{2}{
\nY} \; \ev_j^2\fou{\sol}_j^2 \leq
\frac{2\Opr}{\nY}\;\fou{\sol}_j^2\Opw_j = 2\Opr\zeta\alpha_\nY \leq 1,\]
  
Using the well-known relationship $\rho(\PP_1,\PP_{-1}) \geq 1 - (1/2)KL(\PP_1,\PP_{-1})$ between
the Kullback-Leibler divergence and the Hellinger affinity, we obtain
that $\rho(\PP_1,\PP_{-1})\geq 1/2$.
Using this estimate, \eqref{pr:theo:lower:n:e2} becomes
$\Bigl\{\Ex_{{\theta}}|[\tsol-\sol_\theta]_j|^2+\Ex_{{\theta^{(j)}}}|[\tsol-\sol_{\theta^{(j)}}]_j|^2
\Bigr\} \geq \frac{1}{2}  [\sol]_j^2$, and combining this with~\eqref{eq:11}
implies the result by construction of the solution $\sol$.

\vspace{1em}
\textit{Step (B).} First, we construct two solutions $\sol_\theta\in\solclass$
and operators $\Op_\theta\in \Opclass$ (with $\theta\in\{-1,1\}$) such
that the resulting images $g_\theta$ satisfy $g_{-1}=g_1$. {To
  this end, we
define $\kstars:=\argmax_{j\in\N} \{\hw_j\solw_j^{-1}
\min(1,\nX\Opw_j^{-1})\}$ and $\alpha_\nX:=\zeta
\min(1,\nX^{1/2}\Opw_{\kstars}^{-1/2})$ with
$\zeta:=\min(2^{-1},(1-\Opr^{-1/2}))$.  Observe that
 $1\geq (1-\alpha_\nX)^2\geq
(1-(1-1/d^{1/2}))^2\geq 1/\Opr$ and $1\leq (1+\alpha_\nX)^2\leq
(1+(1-1/\Opr^{1/2}))^2=(2-1/\Opr^{1/2})^2 \leq \Opr$, which
implies $1/\Opr \leq (1+\theta\alpha_\nX)^2\leq
\Opr$.} These inequalities will be used below without further
reference.   We show below
that for each $\theta$ the function {$\sol_\theta:=(1-\theta
\alpha_\nX)  \frac{\solr}{\Opr}\solw_{\kstars}^{-1/2}
\bH_{\kstars}$} belongs to $\solclass$ and that the operator $\Op_\theta$
with the singular values $\ev^\theta_k =
[1+\theta\alpha_\nX\I{k=\kstars}] \;\sqrt{\Opw_k} $ is an element of
$\Opclass$. We obviously have that $\Op_1\sol_f = {(1 - \alpha_\nX^2)(\Opw_{\kstars}/\solw_\kstars)^{1/2}(\solr/\Opr)\bG_\kstars}=\Op_{-1}\sol_{-1}$.

For $\theta\in\{\pm 1\}$, denote by $\PP_\theta$ the joint distribution of
the two sequences $(X_1,X_2,\ldots)$ and $(Y_1,Y_2,\ldots)$, and
let~$\Ex_\theta$ denote the expectation \wrt
$\PP_\theta$.

Applying a reduction scheme as under Step~(A) above, we deduce that for each
estimator $\tsol$ of $\sol$
\begin{align*}
  \sup_{\sol \in \solclass}\sup_{\Op \in \Opclass}
  &\Ex\normV{\tsol -\sol}^2_\hw \geq \max_{\theta\in\{-1,1\}}
  \Ex_{{\theta}}\normV{\tsol-\sol_\theta}^2_\hw\geq\frac{1}{2}
  \Bigl\{\Ex_{1}\normV{\tsol-\sol_1}^2_\hw+\Ex_{-1}\normV{\tsol-\sol_{-1}}^2_\hw
  \Bigr\}.
\end{align*}
Below we show furthermore that 
\begin{equation}\label{pr:theo:lower:m:e1}
  \E_{1}\normV{\tsol-\sol_1}^2_\hw+\E_{-1}\normV{\tsol-\sol_{-1}}^2_\hw \geq {\frac{1}{8}} \normV{\sol_{1}-\sol_{-1}}^2_\hw.
\end{equation}
Moreover, we have $ \normV{\sol_1-\sol_{-1}}_\hw^2 = {4\alpha_\nX^2
(\solr/\Opr)\hw_{\kstars}\solw_{\kstars}^{-1}=
4 \zeta^2(\solr/\Opr)
\hw_{\kstars}\solw_{\kstars}^{-1}
\min\Bigl(1,\frac{\nX}{\edfw_{\kstars}}\Bigr).}$ Combining the last
lower bound with the reduction scheme and the definition of $\kstars$
implies the result of the theorem.

To conclude the proof, it remains to check \eqref{pr:theo:lower:m:e1},
$\sol_\theta\in\solclass$ and $\Op_\theta\in\Opclass$ for both
$\theta$. In order to show $\sol_\theta\in\solclass$,
observe that $\normV{\sol_\theta}^2_\solw=
\solw_{\kstars}|[\sol_\theta]_{\kstars}|^2\leq 
{ \solw_{\kstars}|(1-\theta \alpha_\nX)
(\solr/\Opr) \solw_{\kstars}^{-1/2}|^2 \leq
\solr}$. 

To check that $\Op_\theta\in\Opclass$, it remains to show that
$1/\Opr\leq (\ev^\theta_j)^2/\Opw_j\leq\Opr$ for all $j\geq 1$.  These
inequalities are obviously satisfied for all $j\neq\kstars$, and as
well for $j=\kstars$ by construction of the operator $\Op$. Finally consider
\eqref{pr:theo:lower:m:e1}. As in Step~(A) above by employing the Hellinger affinity
$\rho(\PP_1,\PP_{-1})$ we obtain for any estimator $\tsol$ of $\sol$ that
\begin{equation*}
  \Ex_{1}\normV{\tsol-\sol_{1}}_\hw^2+\Ex_{-1}\normV{\tsol-\sol_{-1}}_\hw^2
 \geq  \frac{1}{2} \normV{\sol_{1}-\sol_{-1}}^2_\hw
  \rho^2(\PP_1,\PP_{-1}). 
\end{equation*}
Next, we bound the Hellinger affinity $\rho(\PP_1,\PP_{-1})$ from below
for all $\nX\in(0,1)$, which proves~\eqref{pr:theo:lower:m:e1}.

Notice that by construction of $\sol_\theta$ and $\Op_\theta$, the
distribution of $X_i$ and $Y_i$ does not depend on $\theta$, except for $X^\theta_\kstars$. It is
thus easily seen that the Kullback-Leibler divergence can be
controlled as follows,
\[KL(\PP_1,\PP_{-1}) ={
\frac{(\ev_\kstars^1 - \ev_\kstars^{-1})^2}{2\nX}  = \frac{2\alpha^2_\nX}{\nX}\,\Opw_\kstars\leq 1}\]
Using $\rho(\PP_1,\PP_{-1}) \geq 1 - (1/2)KL(\PP_1,\PP_{-1})$
again,~\eqref{pr:theo:lower:m:e1} is shown and so is the theorem.
\qed

\subsubsection*{Upper risk bound}
\label{sec:upper-risk-bound}

The following proof uses Lemma~\ref{lem:Rjetc} from the auxiliary
results section~\ref{sec:auxiliary-results} below.

\proofof{Theorem~\ref{theo:upper}} 
Define $\tsol:=\sum_{j=1}^\kstar
[\sol]_j\1\{X_j^2\geq \nX\} \ef_j$ and decompose the risk into
two terms,
\begin{equation}\label{pr:theo:upper:e1}
\Ex\normV{\hsol-\sol}_\hw^2 = 
\Ex\normV{\hsol-\tsol}_\hw^2+\Ex\normV{\tsol-\sol}_\hw^2=: A + B,
\end{equation}
which we bound separately. Consider first $A$ which we decompose
further,
\begin{multline*}
\Ex\normV{\hsol- \tsol}_\hw^2
=  \sum_{j=1}^\kstar \hw_j\E\left[\frac{(Y_j-\E Y_j)^2}{X_j^2} \I{X_j^2\geqslant \nX}\right] \\\hfill
 +\sum_{j=1}^\kstar\hw_j |[\sol]_j|^2
 \E\left[\frac{(X_j-\E X_j)^2}{X_j^2}
 \1\{X_j^2\geqslant \nX\}\right]=: A_1 + A_2.
\end{multline*}
As far as $A_1$ is considered, we use Lemma~\ref{lem:Rjetc}~(iii) from
Section~\ref{sec:auxiliary-results} below and write
\[A_1 = \sum_{j=1}^\kstar \frac{\hw_j\nY}{\E[X_j]^2}\; \E\left[\left(
    \frac{\E[X_j]}{X_j}  \right)^2 \I{X_j^2\geq \nX}  \right]
 \leq 4\Opr\sum_{j=1}^\kstar \frac{\hw_j\nY}{\Opw_j}
 \leq  4\Opr\krate_\nY.\]
As for $A_2$, we apply Lemma~\ref{lem:Rjetc}~(i) and obtain
\[A_2 \leq 8\Opr \sum_{j=1}^\kstar\hw_j |[\sol]_j|^2
                  \min\left(1,\frac{\nX}{\Opw_j}\right)\leq 8d\urate_\nX\]

Consider now $B$ which we decompose further into
\begin{multline*}
\Ex\normV{\tsol- \sol}_\hw^2=\sum_{j\in\N} \hw_j|[\sol]_j|^2 \E[(1-
\1\{1\leq j\leq \kstar\}\1\{X_j^2\geqslant \nX\})^2] \\\hfill
=  \sum_{j > \kstar} \hw_j|[\sol]_j|^2
 +\sum_{j = 1}^\kstar\hw_j |[\sol]_j|^2
 \P\Bigl(X_j^2< \nX\Bigr)=: B_1 + B_2,
\end{multline*}
where $B_1\leq \normV{\sol}^2_\solw
\hw_{\kstar}\solw_{\kstar}^{-1}\leq \solr \krate_\nY$ because
$\sol\in\solclass$.  
Moreover, $B_2\leq 4 \Opr\solr \urate_\nX$ using Lemma~\ref{lem:Rjetc}~(ii).
 The result of the theorem follows now by
combination of the decomposition \eqref{pr:theo:upper:e1} and the
estimates of $A_1,A_2,B_1$ and $B_2$.\qed


%% file: gsm_adaptive_proofs.tex
\subsection{Adaptive estimation (Section \ref{sec:adaptive})}
\label{sec:adapt-estim-proofs}

The proofs in this section use the
Lemmas~\ref{lem:petrala}--~\ref{lem:events} from the auxiliary results
section~\ref{sec:auxiliary-results} below.

\vspace{1ex}

\proofof{Proposition~\ref{prop:concentration}}
Using the model equation $Y_j =
\fou{\im}_j + \sqrt{\nY}\,\xi_j$, we have for all $t\in\cS_k$ that
\begin{equation*}
\fou{\hsol_k - \sol_k}_j
= \frac{\sqrt{\nY}\,\xi_j}{\ev_j} 
+  \left(
 \frac{1}{X_j}\1_{[X_j^2\geq\nX]} -\frac{1}{\ev_j}\right)\sqrt{\nY}\,\xi_j 
+  \left(
 \frac{1}{X_j}\1_{[X_j^2\geq\nX]} -\frac{1}{\ev_j}\right)\fou{\im}_j.
\end{equation*}
Thus, we may decompose the norm $\normV{\hsol_k-\sol_k}_\hw^2$ in
three terms according to 
\begin{align*}
  \normV{\hsol_k - \sol_k}_\hw^2
&\leq 3 \sum_{j=1}^k\frac{\hw_j}{\ev_j} {\nY}\,\xi_j^2
+ 3\sum_{j=1}^k\hw_j \left(
 \frac{1}{X_j}\1_{[X_j^2\geq\nX]} -\frac{1}{\ev_j}\right)^2 {\nY}\,\xi_j^2
+ 3\sum_{j=1}^k\hw_j \left(
 \frac{1}{X_j}\1_{[X_j^2\geq\nX]}
 -\frac{1}{\ev_j}\right)^2\fou{\im}_j^2\\
&=: 3\,\big\{T^{(1)}_k + T^{(2)}_k + T^{(3)}_k\big\}.
\end{align*}
Define the event
\begin{equation*}\label{eq:20}
 \Omega_\nX := \bigg\{ \forall\; 0<j\leq \Ce_\nX^+ \; \bigg|\quad
\Big|\frac{1}{X_j}-
 \frac{1}{\ev_j}\Big| \leq \frac{1}{2\,\ev_j} \quad\text{and}\quad
 X_j^2 \geq \nX\bigg\}.
\end{equation*}
Since  $\I{X_j^2\geq \nX}\I{\Omega_\nX} = \I{\Omega_\nX}$,   it follows
that for all $1\leq j \leq K^+_{\nY,\nX}$ we have
\begin{align*}
  \bigg(\frac{\ev_j}{X_j}\I{X_j^2\geq \nX} - 1\bigg)^2\,\I{\Omega_\nX} =  \ev_j^2\;\I{\Omega_\nX}\,
\bigg|\frac{1}{X_j}-\frac{1}{\ev_j}\bigg|^2
\leq \frac{1}{4}.
\end{align*}
Hence,  
$  T^{(2)}_k\1_{\Omega_\nX} 
\leq \frac{1}{4} T^{(1)}_k $
   for all $1\leq k\leq K_{\nY,\nX}^+$, and thus
\begin{align*}\label{eq:3}
  \begin{split}
    \max_{1\leq k \leq
      K_{\nY,\nX}^+}\left(\normV{\hsol_k-\sol_k}_\hw^2 -
      \frac{1}{6}\pen^a_k\right)_+ \leq 4
    &\sum_{k=1}^{K_{\nY,\nX}^+}\left(\sum_{j=1}^k\frac{\hw_j}{\ev_j}
      {\nY}\,\xi_j^2 - 2\delta_k\nY \right)_+ \\
&+ 3\max_{1\leq k \leq
      K_{\nY,\nX}^+} T^{(2)}_k \1_{\Omega_\nX^c}
+ 3 \max_{1\leq k \leq
      K_{\nY,\nX}^+} T^{(3)}_k.
  \end{split}
\end{align*}
Keeping in mind that $\P[\Omega_\nX^c]\leq C(d)\nX^2$ by virtue of Lemma~\ref{lem:events},
the result follows immediately using Lemmas~\ref{lem:petrala},
\ref{lem:Q}, and \ref{lem:Q2} below.
\qed


\proofof{Proposition~\ref{prop:schachtel}} Let
$\breve\sol_{k}:=\sum_{1\leq j\leq k} \fou{\sol}_j \1\{X_j^2\geq
\nX\}\ef_j$. It is easy to see that
$\normV{\hsol_{k}-\breve\sol_{k}}^2 \leq
\normV{\hsol_{k'}-\breve\sol_{k'}}^2$ for all $ k'\leq k$ and
$\normV{\breve\sol_{k}-\sol}^2\leq \normV{\sol}^2$ for all $k\geq
1$. Thus, using that $1\leq \hk\leq ({\Cy_\nY^\circ}\wedge
\nX^{-1})$, we can write
\begin{align*}
\Ex\normV{\hsol_{\hk}-\sol}_\hw^2\I{\mho_{\nY,\nX}^c} &\leq
2\{\Ex\normV{\hsol_{\hk}-\breve\sol_{\hk}}_\hw^2\I{\mho_{\nY,\nX}^c} +
\Ex\normV{\breve\sol_{\hk}-\sol}_\hw^2\I{\mho_{\nY,\nX}^c}\}\\
&\leq 2\bigg\{ \Ex \normV{\hsol_{({\Cy_\nY^\circ}\wedge \lfloor\nX^{-1}\rfloor)}-\breve\sol_{({\Cy_\nY^\circ}\wedge
     \lfloor\nX^{-1}\rfloor)}}_\omega^2\I{\mho_{\nY,\nX}^c} +
\normV{\sol}_\omega^2\, \P[\mho_{\nY,\nX}^c] \bigg\}.
\end{align*}
 Moreover, using the Cauchy-Schwarz inequality, we
 conclude 
\begin{align*}
  \Ex \normV{\hsol_{({\Cy_\nY^\circ}\wedge  \lfloor\nX^{-1}\rfloor)}&-\breve\sol_{({\Cy_\nY^\circ}\wedge
       \lfloor\nX^{-1}\rfloor)}}_\hw^2\I{\mho_{\nY,\nX}^c} \\
&\leq 2 \nX^{-1}\sum_{1\leq j\leq({\Cy_\nY^\circ}\wedge  \lfloor\nX^{-1}\rfloor)}\omega_j\Bigl\{\Ex(Y_j -\ev_j \fou{\sol}_j)^2\I{\mho_{\nY,\nX}^c} + \Ex(\ev_j \fou{\sol}_j-X_j \fou{\sol}_j)^2\I{\mho_{\nY,\nX}^c}\Bigr\}\\
&\leq 2  \nX^{-1} \Bigl\{ \sum_{1\leq j\leq({\Cy_\nY^\circ}\wedge  \lfloor\nX^{-1}\rfloor)} \hw_j \Bigl[\Ex
  \left(Y_j - \fou{\im}_j\right)^{4}\Bigr]^{1/2}
  \P[\mho_{\nY,\nX}^c]^{1/2}\\
&\hspace{9em}  + \sum_{1\leq <j\leq({\Cy_\nY^\circ}\wedge  \lfloor\nX^{-1}\rfloor)} \hw_j \fou{\sol_{j}}^2[\Ex (X_j-\ev_j)^4]^{1/2}\P[\mho_{\nY,\nX}^c]^{1/2} \Bigr\}  \\
  &\leq 2\sqrt{3}  \nX^{-1} \Bigl\{  ( \nX^{-1}
  {\max_{1\leq j \leq {\Cy_\nY^\circ}}\hw_j}) \nY   +  \nX
  \normV{\sol}_\hw^2\Bigr\} \,\P[\mho_{\nY,\nX}^c]^{1/2},
\end{align*}
which implies
\begin{equation*}
  \Ex\normV{\hsol_{\hk}-\sol}_\hw^2\I{\mho_{\nY,\nX}^c}  \leq
C\bigg\{   \Bigl(  \nX^{-2}
+  \normV{\sol}_\hw^2\Bigr) \,\P[\mho_{\nY,\nX}^c]^{1/2} + \normV{\sol}_\hw^2\, \P[\mho_{\nY,\nX}^c]  \bigg\}.
\end{equation*}
Lemma~\ref{lem:events} below yields, for some 
 $C>0$ depending only on the class $\Opclass$, 
\begin{align*}
  \Ex\normV{\hsol_{\hk}-\sol}_\hw^2\I{\mho_{\nY,\nX}^c} 
\leq C\,\bigg\{ \nX
 + \normV{\sol}_\omega^2{\nX^6}
+ \normV{\sol}_\omega^2{\nX^{12}}  \bigg\}
\end{align*}
which completes the proof due to $\sol\in\solclass$.\qed


%% file: gsm_auxiliary.tex
\subsection{Auxiliary results}
\label{sec:auxiliary-results}

\begin{lem}\label{lem:Rjetc} For every $j\in\N$, 
\ii
  \begin{enumerate}\item 
    $R^{I}_j := \E\bigg[\;\bigg(\frac{\ev_j}{X_j}-1
    \bigg)^2\I{X_j^2\geq \nX} \;\bigg] \leq  \,\min\Big\{1, \frac{8\nX}{\ev_j^2}\Big\}$
  \item $R^{II}_j := \P[X_j^2<
    \nX] \leq \,\min\Big\{1, \frac{4\nX}{\ev_j^2}\Big\}$
  \item $\E\left[\left(\frac{\E[X_j]}{X_j}\right)^2\I{X_j^2\geq\nX}
    \right] \leq 4 $
  \end{enumerate}
\end{lem}
\proof
(i) It is easy to see that
\begin{align}\label{eq:13}
  R^I_j = \E\bigg[ \frac{|X_j - \ev_j|^2}{X_j^2} \;
  \;\I{X_j^2\geq \nX} \bigg]
 \leq \nX^{-1} \var(X_j) =1.
\end{align}
On the other hand, using that $\E[(X_j - \ev_j)^4] = 3\nX^2$, we obtain 
\begin{align*}
 R^I_j  &\leq \E\bigg[ \frac{(X_j-\ev_j)^2}{X_j^2} \;
  \;\I{X_j^2\geq \nX} \; 2
  \bigg\{ \frac{(X_j-\ev_j)^2}{\ev_j^2} +
  \frac{X_j^2}{\ev_j^2} \bigg\} \bigg]\\&\leq
\frac{2\,\E[(X_j-\ev_j)^4]}{\nX\ev_j^2}
+
\frac{2\; \var(X_j)}{\ev_j^2}
=
\frac{8\nX}{\ev^2}.
\end{align*}
Combining with \eqref{eq:13} gives $R^I_j\leq \,\min\Big\{1,
\frac{8\nX}{\ev_j^2}\Big\}$, which completes the proof of (i).\par
(ii) Trivially, $R_j^{II}\leq 1$. 
If $1\leq 4\nX/\ev_j^2$, then obviously $R^{II}_j\leq
\,\min\Big\{1, \frac{4\nX}{\ev_j^2}\Big\}$.
Otherwise, we have $\nX< \ev_j^2/4$ and hence, using \tschebby's
inequality, 
\begin{align*}
  R_j^{II} 
\leq \P[|X_j-\ev_j| > |\ev_j|\,/2\,]
\leq \frac{4\,\var(X_j)}{\ev_j^2}
\leq \,\min\Big\{1, \frac{4\nX}{\ev_j^2}\Big\},
\end{align*}
where we have used that $\var(X_j) = \nX$ for all $j$.\par
(iii)  $\E\left[\left(\frac{\E[X_j]}{X_j}\right)^2\I{X_j^2\geq\nX}
    \right] \leq 2\E\left[ \left(\frac{X_j - \E[X_j]}{X_j}
      \right)^2\I{X_j^2\geq\nX}  + \I{X_j^2\geq\nX}  \right]\leq 4$.
\qed


\begin{lem}\label{lem:NuMu}
Under Assumption~\ref{ass:minreg}, we have  that
 \ii
  \begin{enumerate}\item 
    $  \nY\delta_{\Cy^+_\nY}  \leq 32\,\edfr^2 $ for all $\nY\in(0,1)$,
  \end{enumerate}
and there is a $\nX_0\in(0,1)$ such that for all $\nX <
\nX_0$, we have 
\begin{enumerate}\setcounter{enumi}{1}\item 
  $\min_{1\leq j \leq \Ce^+_\nX} \ev_j^2\geq
  3\nX.$
\end{enumerate}
\end{lem}
\proof
(i) For $\Cy_\nY^+ = 0$, we have $\delta_{\Cy_\nY^+} = 0$ and there is
nothing to show. If $0<\Cy_\nY^+\leq n$, one can show that
$\hw_{\Cy_\nY^+}^+/\edfw_{\Cy_\nY^+}\leq 4\edfr  / (\nY\Cy_\nY^+|\log\nY|)$,
which we use in the following computation:
\begin{align*}
  \delta_{\Cy_\nY^+} &= \Cy_\nY^+ \;
  \frac{\hw_{\Cy_\nY^+}^+}{\edfw_{\Cy_\nY^+}}
 \;\frac{\log((\hw_{\Cy_\nY^+}^+/\edfw_{\Cy_\nY^+}) \vee (\Cy_\nY^+ + 2)
   )}{\log(\Cy_\nY^+ + 2)}
\leq  \frac{4\edfr}{\nY|\log \nY|} 
\; \frac{\log\left( \frac{4\edfr }{\Cy_\nY^+\nY|\log \nY|} \vee (\Cy_\nY^+ + 2)
  \right)}{\log(\Cy_\nY^+ + 2)}\\[1em]
&\leq \nY^{-1}\,
\begin{cases}
  4\edfr & (\log (\nY^{-1}+2)\geq 4\edfr)\\
4\edfr(4\edfr + \log(4\edfr))/(\log (\nY^{-1}+2)) & (\text{otherwise}),
\end{cases}
\end{align*}
which implies $\nY\delta_{\Cy_\nY^+} \leq 4\edfr(4\edfr +
\log(4\edfr))\leq 32\edfr^2$
for all $\nY\in (0, 1)$.
\par
(ii) We have that
\[\min_{1\leq j\leq\Ce_\nX^+} \ev_j^2
\geq \min_{1\leq j\leq\Ce_\nX^+} \frac{\edfw_j}{\edfr} 
\geq \frac{\nX^{1-v_\nX}}{4\edfr^2}\geq 3\nX,\]
where the last step holds for sufficiently small $\nX$ as some
algebra shows.\qed


\begin{lem}\label{lem:petrala}
We have that
\begin{equation*}
  \sum_{ k=1}^{K_{\nY,\nX}^+}\E\bigg( 
  \sum_{j=1}^k \frac{\hw_j}{\ev_j}\nY\xi_j^2 - 2\,\delta^\ev_{k}\nY \bigg)_+ 
\leq 6720\;\nY.
\end{equation*}
\end{lem}

\proof
Representing the expectation of the positive
random variable by the integral over its tail probabilities and using
$\delta_k^\ev\geq\sum_{j=1}^k(\hw_j/\ev_j^2)$, we may write
\begin{multline*}
  \sum_{ k=1}^{K_{\nY,\nX}^+}\E\bigg( 
  \sum_{j=1}^k \frac{\hw_j}{\ev_j}\nY\xi_j^2 - 2\,\delta^\ev_{k}\nY \bigg)_+ 
 \leq 
  \sum_{ k=1}^{K_{\nY,\nX}^+}  
 \int_0^\infty \P\left[\sum_{j=1}^{k}
   \frac{\nY\hw_j}{\ev_j^2}\,(\xi_j^2-1)
       \geq x + 2\nY\delta^\ev_k - \nY\sum_{j=1}^{k}\frac{\hw_j}{\ev_j^2}
     \right] dx
\\  \leq  \sum_{ k=1}^{K_{\nY,\nX}^+}
 \int_0^\infty \P\left[\sum_{j=1}^{k}
   \frac{\nY\hw_j}{\ev_j^2}\,(\xi_j^2-1)
       \geq x + \nY\delta^\ev_k  \right] dx
\end{multline*}
Define $\rho_k := (\nY\hw_k)/\ev_k^2$, $H_k := 4\nY\Delta^\ev_k$, and $B_k := 2\nY^2\sum_{j=1}^{k}\hw_j^2 /
\ev_j^4$. It can be shown (see proof of Proposition~A.1 in \cite{DP:06}) that for all $1\leq k'\leq k$ and $m\geq 2$,
we have
\[\Big|\E\Big[\Big(\frac{\nY\hw_{k'}}{\ev_{k'}^2}(\xi_{k'}^2-1)\Big)^m\Big]\Big|
 \leq m!\,\rho_{k'}^2\, H_k^{m-2}.\]
Hence, the assumption of Theorem~2.8 from \cite{Petrov1995} is
satisfied and splitting up the integral, we get the following bound:
\begin{multline*}
  \sum_{k=1}^{K_{\nY,\nX}^+}\E\bigg( 
  \sum_{j=1}^k \frac{\hw_j}{\ev_j}\nY\xi_j^2 - 2\,\delta^\ev_{k}\nY \bigg)_+ 
\\ \leq  \sum_{k=1}^{K_{\nY,\nX}^+}
\int_0^{B_{k}/H_{k} - \nY\delta^\ev_{k}} \exp\Big( - \frac{(x +
  \nY\delta^\ev_{k})^2}{4B_{k}}\Big) dx
+ \int_{B_{k}/H_{k} - \nY\delta^\ev_{k}}^\infty \exp\Big(-\frac{x +
  \nY\delta^\ev_{k}}{4H_{k}}\Big) dx
\end{multline*}
The second integral is equal to $4H_{k}\exp(-B_{k}/(4H_{k}^2))$. Some computation
shows that the first one is bounded from above by
$4H_{k}\big[\exp\big(-\nY^2(\delta^\ev_{k})^2/(4B_{k})\big) - \exp\big(-B_{k}/(4H_{k}^2)\big)\big] $. Thus, the two
identical terms cancel, and we get 
\begin{align*}
  \sum_{ k=1}^{K_{\nY,\nX}^+}
\E\bigg( 
  \sum_{j=1}^k \frac{\hw_j}{\ev_j}\nY\xi_j^2 - 2\,\delta^\ev_{k}\nY \bigg)_+ 
\leq
  16\;\epsilon\;   \sum_{ k=1}^{K_{\nY,\nX}^+} \Delta^\ev_{k}
  \exp\left(-\frac{(\delta^\ev_{k})^2}{8k(\Delta^\ev_{k})^2}\right).
\end{align*}
To complete the proof, we bound the sum on the right hand side as follows,
\begin{align*}
  \sum_{ k=1}^{K_{\nY,\nX}^+} \Delta^\ev_{k}&
  \exp\left(-\frac{(\delta^\ev_{k})^2}{8k(\Delta^\ev_{k})^2}\right)
\leq  \sum_{k=1}^\infty\exp\Big(-\log(\Delta^\ev_{k}\vee(k+2))
\Big[\frac{k}{8\log(k+2)} -1 \Big]   \Big)\\
&\leq  e\sum_{k=1}^\infty\exp\Big(-
\frac{k}{8\log(k+2)} \Big)
\leq  e\sum_{k=1}^\infty\exp\Big(-
\frac{\sqrt{k}}{8\log(3)} \Big)\\
&\leq  e\int_{0}^\infty\exp\Big(-
\frac{\sqrt{x}}{8\log(3)} \Big)dx=128\log^2(3)\,e,
\end{align*}
where we have used 
$\log(k+2)\leq \log(3)\sqrt{k}$ for all $k\geq 1$. 
\qed


\begin{lem}\label{lem:Q}
For every $k\in\N$ and $\nX\in(0,1)$,
\[\E\bigg[\sum_{j=1}^k\hw_j\fou{\im}_j^2\left(\frac{1}{X_j}\1_{[X_j\geq\nX]}
- \frac{1}{\ev_j}\right)^2   \bigg] \leq 
8\;\Opr\;\solr\;\urate_\nX(\solw,\Opw,\hw).\]
\end{lem}
\proof Firstly, as $\sol\in\solclass$, it is easily seen that
\begin{align*}
\E\bigg[\sum_{j=1}^k\hw_j\fou{\im}_j^2\left(\frac{1}{X_j}\1_{[X_j\geq\nX]}
- \frac{1}{\ev_j}\right)^2   \bigg]
\leq
\solr\;\max_{1\leq j\leq k}\;\frac{\hw_j}{\solw_j}\,\E[|R_j|^2],
\end{align*}
where $R_j$ is defined as 
\begin{equation}
R_j := \left(\frac{\ev_j}{X_j}\I{X_j^2\geq\nX^2} -1\right).\label{eq:6}
\end{equation}
 In view of the definition of
$\urate_\nX$ in Theorem~\ref{theo:lower}, the result follows from $\E[|R_j|^2] \leq \Opr\,\min\Big\{1,
\frac{8\nX}{\Opw_j}\Big\}$, which is a consequence of the decomposition
\begin{align}\label{eq:18}
  \begin{split}
    \E|R_j|^2  
    = \E\bigg[\;\bigg(\frac{\ev_j}{X_j}-1
    \bigg)^2\I{X_j^2\geq \nX} \;\bigg]
    +\P[X_j^2< \nX] 
  \end{split}
\end{align}
and Lemma~\ref{lem:Rjetc}.
\qed


\begin{lem}\label{lem:Q2} We have that
  \[\E\bigg[\sum_{j=1}^{K_{\nY,\nX}^+}
  \hw_j\left(\frac{1}{X_j}\1_{[X_j\geq\nX]}-
   \frac{1}{\ev_j}\right)^2{\nY}\xi_j^2 \1_{\Omega_\nX^c}\bigg]
{\leq 64\,\Opr^3 (\P[\Omega_\nX^c])^{1/2}}
.\]
\end{lem}
\proof  
Given $R_j$ from~\eqref{eq:6},  we begin our proof observing that
\begin{align*}
\E\bigg[\sum_{j=1}^{K_{\nY,\nX}^+}
  \hw_j\left(\frac{1}{X_j}\1_{[X_j\geq\nX]}-
   \frac{1}{\ev_j}\right)^2\sqrt{\nY}\xi_j^2 \1_{\Omega_\nX^c}\bigg]
  &\leq \nY \sum_{j=1}^{K_{\nY,\nX}^+}\frac{\hw_j}{\ev_j^2}
  \;\E[|R_j|^2\1_{\Omega_\nX^c}],
\intertext{{where we have used the independence of  $X$ and $Y$ and
  $\var(Y_j) = \nY$. Since $\Opr\delta_k^\Opw \geq \sum_{j=1}^k \frac{\omega_j}{\ev_j^2}$ for all $\Op\in\Opclass$, the Cauchy-Schwarz inequality yields}}
  \E\bigg[\sum_{j=1}^{K_{\nY,\nX}^+}
  \hw_j\left(\frac{1}{X_j}\1_{[X_j\geq\nX]}-
   \frac{1}{\ev_j}\right)^2{\nY}\xi_j^2 \1_{\Omega_\nX^c}\bigg]
&{\leq\; d\,(\P[\Omega_\nX^c])^{1/2}\;{\nY\delta^\Opw_{N^+_\nY}}
\;\max_{0<j\leq N^+_\nY}(\E[|R_j|^4])^{1/2}.}
\end{align*} 
Proceeding analogously to~(\ref{eq:13}) and~(\ref{eq:18}), one can
show that $\E[|R_j|^4]\leq 4$. The result follows then using the
  definition of $\Cy^+_\nY$.\qed


\begin{lem}\label{lem:events}
For $k\in\N$,  define the events
\[\tOmega_k := \bigg\{\Big|\frac{X_j}{\ev_j} - 1\Big| \leq \frac{1}{3} \quad
\forall\, 1\leq j \leq k   \bigg\} \]
and suppose that  Assumption~\ref{ass:minreg} holds. 
For all $\nY,\nX\in(0,1)$ , we have 
\begin{enumerate}\ii
\item $ \Omega_\nX \subseteq \{\pen^+_k\leq \hpen_k \leq 30\pen^+_k \quad
\forall\;1\leq k \leq K^+_{\nY,\nX}\}$,
\item $\tOmega_{\Ce_\nX^++1}\subseteq \{K^-_{\nY,\nX} \leq
\hK_{\nY,\nX} \leq K^+_{\nY,\nX}\}$,
\item $\P[\tOmega_{\Ce_\nX^+}^c]\leq C(\Opr)\,\nX^2$ \;and\;
  $\P[\Omega_\nX^c]\leq C(\Opr)\,\nX^2$.
\end{enumerate}
If additionally condition~\eqref{eq:Mplusone} holds, then
\begin{enumerate}\ii\setcounter{enumi}{3}
\item  $\P[\mho_{\nY,\nX}^c] 
\leq  C(\Opw,\Opr) \nX^6$.
\end{enumerate}
\end{lem}
\proof Consider (i).  Notice first that
$  \delta^\ev_k \leq 
 \delta_k^\Opw\, \Opr\,\zeta_\Opr$ for all $k\geq
    1$
with $\zeta_\Opr:= (\log (3d))/(\log 3)$. 
Observe that on $\Omega_\nX$ we have  $ (1/2)\Delta^\ev_k\leq \Delta^X_k\leq(3/2)\Delta^\ev_k$ for all $1\leq k\leq
\tCe_\nX$
and hence $(1/2)[\Delta^\ev_k\vee(k+2)]\leq[ \Delta^X_k\vee
(k+2)]\leq(3/2)[\Delta^\ev_k\vee(k+2)]$, which implies
\begin{align*}
  \begin{split}
    (1&/2) k \Delta^\ev_k\Bigl(\frac{ \log[ \Delta^\ev_k\vee
      (k+2)]}{\log(k+2)}\Bigr)\Bigl(1-\frac{\log 2}{\log
      (k+2)}\frac{\log (k+2)}{\log (\Delta^\ev_k \vee [k+2])}\Bigr) 
\\[1ex]&\leq
    \delta^X_k \leq (3/2) k \Delta^\ev_k\Bigl(\frac{\log (\Delta^\ev_k
      \vee [k+2])}{\log( k+2)}\Bigr) \Bigl(1+\frac{\log3/2}{\log(k+2)
    } \frac{\log (k+2)}{\log (\Delta^\ev_k \vee [k+2])}\Bigr).
  \end{split}
\end{align*}
Using $ {\log (\Delta^\ev_k
  \vee (k+2))}/{\log( k+2)}\geq 1$, we conclude from the last estimate that
\begin{align*}
  \begin{split}
{ \delta^\ev_k/10\leq}
(\log 3/2)/(2 \log 3) \delta^\ev_k&\leq (1/2) \delta^\ev_k[1-(\log
2)/\log(k+2)] \leq {\delta^X_k} \\ &\leq (3/2) \delta^\ev_k[1+
(\log3/2)/\log (k+2)]\leq { 3\delta^\ev_k}.
\end{split}
\end{align*} 
It follows that on $\Omega_\nX$ we have $\pen^+_k\leq \hpen_k\leq 30
\pen^+_k$  for all $1\leq k\leq \Ce_\nX^+$ as desired.\\[1ex]
Proof of (ii). 
Denoting by $X$ the random sequence $(X_j)_{j\geq1}$, define sequences $N_\nY^- :=
  N_\nY^{\sqrt{\Opw/(4\Opr)}}$, $M_\nX^- :=
  M_\nX^{\sqrt{\Opw/(4\Opr)}}$ and   $\hN_\nY := N_\nY^{X}$, $\hM_\nX
  := M_\nX^X$. Note that by definition, $K_{\nY,\nX}^- = N_\nY^- \wedge M_\nX^-$ and
$\hK_{\nY,\nX} = \hN_\nY \wedge \hM_\nX$.
Define further the events  $\Omega_I:=\{ K^-_{\nY,\nX}> \hK_{\nY,\nX}\}$ and
 $\Omega_{II}:=\{ \hK_{\nY,\nX}>K^+_{\nY,\nX}\}$. 
Then we have $\{K^-_{\nY,\nX} \leq
\hK_{\nY,\nX} \leq K^+_{\nY,\nX}\}^c=\Omega_I\cup \Omega_{II}$.  
Consider $\Omega_I = \{ \hCy_\nY< K^-_{\nY,\nX}\} \cup \{
\hCe_\nX< K^-_{\nY,\nX}\}$ first.  By definition of $\Cy^-_\nY$,
we have that $\min_{1\leq j\leq\Cy^-_\nY}
\frac{\ev_j^{2}}{j\,\hw_j^+}\geq 4 \,\nY |\log \nY|$, which implies,
keeping in mind that $K^-_{\nY,\nX}\leq \Cy^-_{\nY,\nX}$,
   \begin{multline*}
     \{ \hCy_\nY< K^-_{\nY,\nX}\}\subset \bigg\{ \exists\, 1\leq
     j\leq K^-_{\nY,\nX}\,\bigg|\,
     \frac{X_j^2}{j\,\hw_j^+}<\nY|\log\nY|\bigg\}\\
\subset
     \bigcup_{1\leq j\leq K^-_{\nY,\nX}}\bigg\{
     \frac{X_j}{\ev_j}\leq \frac{1}{2}\bigg\}
\subset
     \bigcup_{1\leq j\leq K^-_{\nY,\nX}}\bigg\{
     \left|\frac{X_j}{\ev_j}-1\right|\geq \frac{1}{2}\bigg\}.
  \end{multline*}
One can see that from $\min_{1\leq j\leq\Ce^-_\nX}
\ev_j^{2}\geq 4 \nX^{1- v_\nX}$ it follows in the same way that
   \begin{equation*}
     \Big\{ \hCe_\nX < K^-_{\nY,\nX}\Big\}\subset 
     \bigcup_{1\leq j\leq K^-_{\nY,\nX}}\bigg\{ \left|\frac{X_j}{\ev_j}-1\right|\geq \frac{1}{2}\bigg\}.\hfill
  \end{equation*} 
Therefore, $\Omega_I \subseteq \bigcup_{1\leq j\leq 
  \Ce^+_\nX}\Bigl\{ |X_j/\ev_j-1|\geq 1/2\Bigr\}\subseteq \tOmega_{\Ce_\nX^++1}^c$, since
$\Ce_\nX^-\leq \Ce_\nX^+$. 

Consider $\Omega_{II} =  \{ \hCy_\nY> K^+_{\nY,\nX}\} \cap  \{ \hCe_\nX>
 K^+_{\nY,\nX}\}$.
 In case $K^+_{\nY,\nX}= \Cy_\nY^+$, note that by definition of
 $\Cy_\nY^+$, we have $\nY|\log \nY| / 4 \geq
\frac{\ev_{\Cy_\nY^+ +1}^2}{(\Cy_\nY^+ +1)\,\hw_{\Cy_\nY^+ +1}^+}$, such that 
 \begin{align*}
   \Omega_{II}\subseteq\{\hCy_\nY> \Cy_\nY^+\}&\subset
 \Bigl\{ \forall 1\leq j\leq
   \Cy_\nY^+ + 1\quad\bigg|\quad \frac{X_j^2}{j\,\hw_j^+}\geq \nY|\log\nY|\Bigr\}\\[1ex] &\subset
   \Biggl\{ \frac{X_{\Cy_\nY^++1}}{\ev_{\Cy_\nY^++1}}\geq
   2\Biggr\}\subset \Biggl\{ \bigg|\frac{X_{\Cy_\nY^++1}}{\ev_{\Cy_\nY^++1}}-1\bigg|\geq
   1\Biggr\}.
  \end{align*}
 In case $K^+_{\nY,\nX}= \Ce_\nX^+$, it follows analogously from $\nX^{1-v_\nX} \geq
 4 \max_{j\geq \Ce_\nX^++1} \ev_j^2$ that
   \begin{equation*}
     \Omega_{II}\subset\{ \hCe_\nX> \Ce_\nX^+\}
     \subset
     \Bigl\{| X_{\Ce_\nX^++1}/\ev_{\Ce_\nX^++1}-1|\geq 1\Bigr\}.\hfill
  \end{equation*} 
  Therefore, we have $\Omega_{II} \subseteq \Bigl\{ |
  X_{K_{\nY,\nX}^++1}/\ev_{K_{\nY,\nX}^++1}-1 | \geq
  1\Bigr\}\subseteq \tOmega_{\Ce^+_{\nX}+1}^c$ and (ii) is shown. \\[1ex]

Proof of (iii).    For
  $Z\sim\cN(0,1)$ and $z\geq 0$, one has $\P[Z>z]\leq (2\pi
  z^2)^{-1/2}\exp(-z^2/2)$. Hence, there is a constant
  $C(\Opr)$ depending on $\Opr$ such that for every $ 1\leq j
  \leq\Ce^+_\nX$, 
\begin{equation*}
 \label{eq:7}
\P[|X_j/\ev_j -1
|>1/3]
\leq
C(\Opr)\,\left(\frac{\nX}{{\Opw_{\Ce_{\nX}^+}}}\right)^{1/2}\exp\bigg(-
\frac{\Opw_{\Ce_\nX^+}}{18\nX\Opr} \bigg).
 \end{equation*}
Consequently, as $\Ce_\nX^+\leq\nX^{-1}$ and $\Opw_{\Ce_\nX^+}>
\nX^{1-v_\nX} / (4\Opr)$, we have  
\[\P[\tOmega_{\Ce_\nX^+}^c]
\leq C(\Opr) \nX^{2-v_\nX}\exp\bigg(- \frac{\nX^{-v_\nX}}{72\Opr^2} \bigg) \]
 which implies $\P[\tOmega_{\Ce_\nX^+}^c]\leq C(\Opr)\,\nX^2$ using
 that $\nX^{v_\nX}\,|\log\nX|\to 0 $ as $\nX\to 0$. 
As for the second assertion in~(iii), 
we distinguish the cases $\nX\leq\nX_0$ and
  $\nX>\nX_0$, where $\nX_0$ is the constant from
  Lemma~\ref{lem:NuMu}~(ii) depending only on $\Opr$. The assertion is trivial
  for $\nX>\nX_0$ (keeping in mind that $\P[\Omega_{\nX}^c]\leq
  \nX_0^{-2}\nX^2$). Consider the case $\nX\leq\nX_0$, where
  $\ev_j^2\geq 3\nX$ for all $1\leq j \leq \Ce_\nX^+$ due to Lemma~\ref{lem:NuMu}~(ii).
 This yields for the complement of $\Omega_\nX$
 \[\Omega_\nX^c = \bigg\{\exists\;  1\leq j\leq  \Ce_\nX^+  \quad \bigg|\quad
 \Big|\frac{\ev_j}{X_j}-1\Big| > \frac{1}{2} \quad\text{or}\quad
 X_j^2 < \nX \bigg\} \subseteq 
\bigg\{\exists\; 1\leq j\leq  \Ce_\nX^+ \;\bigg|\;
\bigg|\frac{X_j}{\ev_j} -1
\bigg|> \frac{1}{3}\bigg\}=\tOmega^c_{\Ce_\nX^+}.\] 
It follows
with assertion (ii) that $\mho_{\nY,\nX}^c \subseteq
\tOmega_{\Ce_\nX^+}^c$ for all $\nX\leq\nX_0$, implying the second
assertion of (iii).
\medskip

Proof of (iv). Following the proof of (iii) and using that
$\Ce_\nX^++1\leq\nX^{-1}$, we obtain
\begin{equation}
\P[\tOmega_{\Ce_\nX^++1}^c]
\leq C(\Opr) (\nX\Opw_{\Ce_{\nX}^++1})^{-1/2}\exp\bigg(-
\frac{\Opw_{\Ce_\nX^++1}}{18\nX\Opr} \bigg).\label{eq:1}
\end{equation}
Note that
$\tOmega_{\Ce_\nX^++1} \subseteq \Omega_\nX$, since trivially
$\tOmega_{\Ce_\nX^++1} \subseteq \tOmega_{\Ce_\nX^+}$. Thus,~\eqref{eq:1}  implies assertion~(iv)  by virtue of condition~\eqref{eq:Mplusone}.\qed


%% file: gsm.bbl
\begin{thebibliography}{}

\bibitem[Barron et~al., 1999]{BBM:99}
Barron, A., Birg\'e, L., and Massart, P. (1999).
\newblock {Risk bounds for model selection via penalization.}
\newblock {\em Probability Theory and Related Fields}, 113:301--413.

\bibitem[Cavalier et~al., 2002]{CGP:02}
Cavalier, L., Golubev, G., Picard, D., and Tsybakov, A. (2002).
\newblock {Oracle inequalities for inverse problems.}
\newblock {\em Ann. Stat.}, 30:843--874.

\bibitem[Cavalier and Hengartner, 2005]{CH:05}
Cavalier, L. and Hengartner, N.~W. (2005).
\newblock {Adaptive estimation for inverse problems with noisy operators.}
\newblock {\em Inverse Problems}, 21:1345--1361.

\bibitem[Dahlhaus and Polonik, 2006]{DP:06}
Dahlhaus, R. and Polonik, W. (2006).
\newblock {Nonparametric quasi-maximum likelihood estimation for Gaussian
  locally stationary processes.}
\newblock {\em Ann. Stat.}, 34:2790--2824.

\bibitem[Efromovich, 1997]{Efr:97}
Efromovich, S. (1997).
\newblock {Density estimation for the case of supersmooth measurement error.}
\newblock {\em Journal of the American Statistical Association}, 92:526--535.

\bibitem[Ermakov, 1990]{0729.62089}
Ermakov, M. (1990).
\newblock {On optimal solutions of the deconvolution problem.}
\newblock {\em Inverse Probl.}, 6(5):863--872.

\bibitem[Fan, 1991]{Fan:91}
Fan, J. (1991).
\newblock {On the optimal rates of convergence for nonparametric deconvolution
  problems.}
\newblock {\em The Annals of Statistics}, 19:1257--1272.

\bibitem[Goldenshluger and Lepski, 2011]{GL:10}
Goldenshluger, A. and Lepski, O. (2011).
\newblock {Bandwidth selection in kernel density estimation: oracle
  inequalities and adaptive minimax optimality.}
\newblock {\em Ann. Stat.}, 39(3):1608--1632.

\bibitem[Hoffmann and Reiss, 2008]{HR:08}
Hoffmann, M. and Reiss, M. (2008).
\newblock {Nonlinear estimation for linear inverse problems with error in the
  operator.}
\newblock {\em The Annals of Statistics}, 36:310--336.

\bibitem[Johnstone and Silverman, 1990]{JS:90}
Johnstone, I.~M. and Silverman, B.~W. (1990).
\newblock {Speed of estimation in positron emission tomography and related
  inverse problems.}
\newblock {\em Ann. Stat.}, 18(1):251--280.

\bibitem[Mair and Ruymgaart, 1996]{MairRuymgaart1996}
Mair, B.~A. and Ruymgaart, F.~H. (1996).
\newblock {Statistical inverse estimation in Hilbert scales.}
\newblock {\em SIAM Journal on Applied Mathematics}, 56(5):1424--1444.

\bibitem[Math\'e and Pereverzev, 2001]{MP:01}
Math\'e, P. and Pereverzev, S.~V. (2001).
\newblock {Optimal discretization of inverse problems in Hilbert scales.
  Regularization and self-regulari\-za\-tion of projection methods.}
\newblock {\em SIAM J. Numer. Anal.}, 38(6):1999--2021.

\bibitem[Neumann, 1997]{Neu:97}
Neumann, M.~H. (1997).
\newblock {On the effect of estimating the error density in nonparametric
  deconvolution.}
\newblock {\em Journal of Nonparametric Statistics}, 7:307--330.

\bibitem[Petrov, 1995]{Petrov1995}
Petrov, V.~V. (1995).
\newblock {\em Limit theorems of probability theory. Sequences of independent
  random variables.}
\newblock Oxford Studies in Probability. Clarendon Press., Oxford, 4. edition.

\bibitem[Stefanski and Carroll, 1990]{SC:90}
Stefanski, L. and Carroll, R.~J. (1990).
\newblock {Deconvoluting kernel density estimators.}
\newblock {\em Statistics}, 21:169--184.

\end{thebibliography}
